\documentclass[12pt]{article}
\usepackage{amssymb}
\usepackage{latexsym}
\def\part#1{\frac{\partial\phantom{q}}{\partial#1}}

\newenvironment{rmk}{\begin{trivlist}\item[]{\bf Remark:} }
{\end{trivlist}}
\newenvironment{ex}{\begin{trivlist}\item[]{\bf Example:} }
{\end{trivlist}}
\newenvironment{prf}{\begin{trivlist}\item[]{\bf Proof:} }
{\hfill $\Box$ \end{trivlist}}
\newenvironment{lemprf}{\begin{trivlist}\item[]{\bf Proof:} }
 {\end{trivlist}}
\newtheorem{thm}{Theorem}
\newtheorem{definition}{Definition}
\newtheorem{prp}[thm]{Proposition}
\newtheorem{lem}[thm]{Lemma}

\newcommand{\lie}[1]{\mathfrak{#1}}
\def\End{\mathop{\rm End}\nolimits}

\def\tr{\mathop{\rm tr}\nolimits}

\def\Diff{\mathop{\rm Diff}\nolimits}

\newcommand{\R}{\mathbf{R}}
\newcommand{\C}{\mathbf{C}}
\newcommand{\K}{\mathbf{H}}

\newcommand{\CP}{{\mathbf C}{\rm P}}
\newcommand{\RP}{{\mathbf R}{\rm P}}
\begin{document}
\title{The geometry of three-forms in six and seven dimensions}
 \author{Nigel Hitchin\\[5pt]
\itshape  Mathematical Institute\\
\itshape 24-29 St Giles\\
\itshape Oxford OX1 3LB\\
\itshape UK\\
 hitchin@maths.ox.ac.uk}
\maketitle
\begin{abstract}
\noindent We study the special algebraic properties of alternating $3$-forms in $6$ and $7$ dimensions and  introduce a diffeomorphism-invariant functional on the space of differential 
$3$-forms on a closed  manifold $M$ in these dimensions. Restricting the functional to a de Rham cohomology class in $H^3(M,\R)$, we find that a critical point which is generic in a suitable sense defines in the $6$-dimensional case  a complex threefold with trivial canonical bundle and in $7$ dimensions a Riemannian manifold with holonomy $G_2$. This approach gives a direct method of showing that an open set in $H^3(M,\R)$, with a special geometry, is a local moduli space for these structures.
 \end{abstract}

\section{Introduction}
This paper arose from the author's interest in the geometry of $3$-forms on a manifold. Exterior differential forms of degree $2$ are much studied, in particular a symplectic manifold is defined by a closed non-degenerate $2$-form. Our starting point here is the fact that in dimensions $6$ and $7$ the notion of non-degeneracy for {\it three-forms} makes sense too. By this we mean that if $W$ is a real vector space of dimension $6$ or $7$, then the group $GL(W)$ has an open orbit, a fact that has been known for nearly a century \cite{R}. The orbits we are particularly interested in have stabilizer conjugate to $SL(3,\C)$ when $\dim W=6$, and the exceptional Lie group $G_2$ when $\dim W=7$. In this way a single  $3$-form $\Omega$ on a manifold $M$ defines a reduction of the structure group of $M$ to one of these groups. When we look closely at the mechanism for this reduction, there appears naturally a volume form algebraically defined by any $3$-form. Integrating this form gives a diffeomorphism-invariant functional $\Phi$ on the space of $3$-forms, and this is our main object of study.

Following the analogy with Hodge theory, we restrict this functional to {\it closed} forms on $M$ in a given de Rham cohomology class and look for critical points. What we find is that if the critical point $\Omega$ is a $3$-form which lies in the open orbit everywhere, then the reduction of structure group is integrable in the sense that for $\dim M=6$ we obtain a complex threefold with trivial canonical bundle, and for $\dim M=7$ a Riemannian manifold with holonomy in $G_2$. These  geometrical structures are therefore simply critical points of the invariant functionals.

Because the functional is diffeomorphism-invariant, every critical point lies on an orbit of critical points and so $\Phi$ is never a Morse function. However, we show that formally it is a Morse-Bott function -- its  Hessian is non-degenerate transverse to the orbits of $\Diff(M)$. In the $6$-dimensional case this requires the assumption that the complex threefold satisfies the $\partial\bar\partial$-lemma, and in the $G_2$ case we need an analogue of this lemma, which is essentially contained in \cite{Joyce}. This nondegeneracy can be used, together with a standard use of the Banach space implicit function theorem, to give a direct and easy proof that the moduli space of complex structures together with holomorphic $3$-forms on a $6$-manifold is locally an open set in $H^3(M,\R)$. In particular we see that the moduli space of complex structures is unobstructed. The novelty of our approach is that the flat structure on this moduli space is apparent from the very beginning, and the complex structure is defined in a secondary manner. This is the opposite point of view from the conventional use of Kodaira-Spencer theory as in the work of Tian and Todorov \cite{Tian1},\cite{Tod}. For us the complex structure on the moduli space  carries with it the natural special pseudo-K\"ahler structure whose existence is an important ingredient in mirror symmetry \cite{Cand}. We carry out the same process for $G_2$-manifolds and show that the moduli space is isomorphic to an open set in $H^3(M,\R)$ with an indefinite metric of Hessian type, completing a line of thought due to Joyce \cite{Joyce}.

The structure of the paper is as follows. In Section 2 we consider the linear algebra of the vector space $\Lambda^3W^*$ where $W$ is $6$-dimensional. The essential point is that over the complex numbers, a generic $3$-form is the sum of two decomposable ones. In Section 3, we see this algebra from the point of view of symplectic geometry, regarding $\Lambda^3W^*$ as a symplectic vector space under the action of $SL(W)$. This viewpoint is extremely useful for studying the variational problem. Section 4 is a detour into the realm of self-duality. With an inner product of signature $(5,1)$ on $W$, we can define self-dual and anti-self-dual $3$-forms, and some of the linear algebra assumes a much more concrete form. Moreover, we see a setting here for the equation of motion for a self-interacting self-dual tensor, a nonlinear equation of some current interest to physicists. In Section 5 we introduce the invariant functional on $3$-forms on a $6$-manifold, and relate the critical points to integrable complex structures. In Section 6 we prove Morse-Bott nondegeneracy formally, and  then use a Sobolev space model to prove rigorously that an open set in $H^3(M,\R)$ is a local moduli space. Section 7 is concerned with carrying out the parallel programme for $G_2$-manifolds.

The author wishes to thank Robert Bryant for explaining the algebra behind $SL(3,\C)$ and $G_2$ and Patrick Baier for  useful discussions.

\section{Linear algebra}

\subsection{The complex case}

Let $V$ be a $6$-dimensional complex vector space and $\Lambda^3V^*$ the $20$-dimensional
 vector space  of alternating multilinear $3$-forms on $V$. Take $\Omega\in \Lambda^3V^*$ and $v\in V$ and the interior product $\iota(v)\Omega\in  \Lambda^2V^*$. Then
$\iota(v)\Omega\wedge \Omega\in \Lambda^5V^*.$
On the other hand, the natural exterior product pairing
$V^*\otimes \Lambda^5V^*\rightarrow \Lambda^6V^*$
provides an isomorphism
$$A:\Lambda^5V^*\cong V\otimes \Lambda^6V^*$$
and using this we  define a linear transformation
$K_{\Omega}:V\rightarrow V\otimes \Lambda^6V^*$
by
\begin{equation}
K_{\Omega}(v)=A(\iota(v)\Omega\wedge \Omega)
\label{Komega}
\end{equation}
\begin{definition}
Define $\lambda(\Omega)\in (\Lambda^6V^*)^2$ by 
$$\lambda(\Omega)=\frac{1}{6}\tr K_{\Omega}^2$$
\end{definition}
Using $\lambda(\Omega)$ we have the following characterization of ``non-degenerate"
elements in $\Lambda^3V^*$: \begin{prp}
\label{decomprop}
 For $\Omega\in \Lambda^3V^*$, $\lambda(\Omega)\ne 0$ if and only if $\Omega=\alpha+\beta$ where $\alpha,\beta$ are decomposable  and $\alpha\wedge \beta\ne 0$. The $3$-forms $\alpha,\beta$ are unique up to ordering.
\end{prp}
\begin{prf} Let $v_1,\dots,v_6$ be a basis of $V$ and $\theta_1,\dots,\theta_6\in V^*$ the dual basis. Take 
$$\varphi=\theta_1\wedge \theta_2\wedge \theta_3+\theta_4\wedge \theta_5\wedge \theta_6.$$
Let $\epsilon=\theta_1\wedge \dots \wedge \theta_6$ be the associated basis vector for $\Lambda^6V^*$. We find easily that
\begin{equation}
K_{\varphi}v_i=v_i\epsilon\quad (i=1,2,3),\qquad K_{\varphi}v_i=-v_i\epsilon\quad (i=4,5,6),
\label{K}
\end{equation}
and 
\begin{equation}
\lambda(\varphi)=\epsilon^2
\label{L}
\end{equation} 
\vskip .25cm
Now if $\Omega=\alpha + \beta$ where $\alpha,\beta$ are decomposable 
then $\alpha=\xi_1\wedge \xi_2\wedge \xi_3$ and $\beta=\eta_1\wedge \eta_2\wedge \eta_3$ and 
the condition $\alpha\wedge \beta\ne 0$ implies
 that $\{\xi_1, \xi_2,\xi_3,\eta_1,\eta_2,\eta_3\}$ forms a basis for $V^*$.
  Thus by an element of $GL(V)$, $\Omega$ can be transformed to
  $\varphi$.
  Since  $\lambda(\varphi)$ is non-zero from (\ref{L}), so is 
  $\lambda(\Omega)$. \vskip .25cm
If we transform $\varphi$ by $K_{\varphi}$  we find
$$K^*_{\varphi}\varphi=(\theta_1\wedge \theta_2\wedge \theta_3-\theta_4\wedge \theta_5\wedge \theta_6)\epsilon^3$$
and so 
$$K^*_{\varphi}\varphi+\epsilon^3 \varphi=2(\theta_1\wedge\theta_2\wedge\theta_3)\epsilon^3$$
Thus we see that $K_{\varphi}^*\varphi+{\lambda(\varphi)}^{3/2}\varphi\in \Lambda^3V^*\otimes (\Lambda^6V^*)^3$ is decomposable, for each choice of square root of $\lambda(\varphi)$. Moreover, we have
\begin{equation}
K^*_{\varphi}\varphi \wedge \varphi=2\epsilon^4=2\lambda(\varphi)^2
\label{KOO}
\end{equation}
Let $G\subset GL(V)$ be the stabilizer of $\varphi$. Each element of $G$ 
 commutes  with $K_{\varphi}$ and so preserves or interchanges the 
 two subspaces spanned by $v_1,v_2,v_3$ and $v_4,v_5,v_6$ respectively. The identity 
 component $G_0$ preserves them, and the nonvanishing $3$-forms
  $\theta_1\wedge \theta_2\wedge \theta_3$ and $\theta_4\wedge \theta_5\wedge \theta_6$
   defined on them, and so is isomorphic to a subgroup of $SL(3,\C)\times SL(3,\C)$. Thus
$$\dim G_0 \le \dim (SL(3,\C)\times SL(3,\C))=16.$$
But $\dim GL(V)=36$ so the dimension of the orbit is at least $36-16=20$. Since 
$\dim \Lambda^3V^*=20$, the orbit is open and the stabilizer must actually be equal to 
$SL(3,\C)\times SL(3,\C)$. 

For the converse, note that the algebraic condition (\ref{KOO}) and
the decomposability of $K_{\Omega}^*\Omega+{\lambda(\Omega)}^{3/2}\Omega$ hold 
on the open orbit and therefore hold everywhere. Thus in general  we have
decomposable (possibly zero) forms $\alpha,\beta \in \Lambda^3V^*$  with \begin{eqnarray*}
{\lambda(\Omega)}^{3/2}\Omega+K_{\Omega}^*\Omega&=&2{\lambda(\Omega)}^{3/2}\alpha\\
{\lambda(\Omega)}^{3/2}\Omega- K_{\Omega}^*\Omega&=&2{\lambda(\Omega)}^{3/2}\beta
\end{eqnarray*}
and so if $\lambda(\Omega)\ne 0$
\begin{equation}
\Omega=\alpha+\beta,\qquad K_{\Omega}^*\Omega=\lambda(\Omega)^{3/2}(\alpha - \beta)
\label{AB}
\end{equation}
 From (\ref{KOO}), $K^*_{\Omega}\Omega \wedge \Omega=2\lambda(\Omega)^2\ne 0$ but also from (\ref{AB}) 
$$K^*_{\Omega}\Omega \wedge \Omega=2\lambda(\Omega)^{3/2}\alpha\wedge\beta$$ so $\alpha \wedge \beta \ne 0$ as required.
\vskip .25cm
By construction, $\alpha$ and $\beta$ are unique given the choice of square root of $\lambda(\Omega)$.
\end{prf}
\begin{rmk}
 The proposition tells us that the open set $\lambda(\Omega)\ne 0$ is the orbit of the $3$-form $\varphi=\theta_1\wedge \theta_2\wedge \theta_3+\theta_4\wedge \theta_5\wedge \theta_6$ under the action of $GL(V)$. We can therefore deduce properties of $\Omega$ from those of $\varphi$. For example, from (\ref{K}) and (\ref{L}) it follows that
\begin{eqnarray}
\tr K_{\Omega}&=&0\label{trace}\\
K^2_{\Omega}&=&\lambda(\Omega)1\label{Ksquare}
\end{eqnarray}
\end{rmk}

\subsection{The real case}

Now suppose that $W$ is a {\it real} $6$-dimensional vector space and $\Omega \in \Lambda^3W^*$. 
In this case $\lambda(\Omega)\in (\Lambda^6W^*)^2$ is real. If $L$ is a real one-dimensional 
vector space, we say that a vector $u\in L\otimes L=L^2$,  is positive ($u>0$) if 
$u=s\otimes s$ for some $s\in L$ and negative  if $-u>0$. 
\begin{prp}
\label{signprop}
  Suppose that $\lambda(\Omega)\ne 0$ for $\Omega\in \Lambda^3W^*$. Then
\begin{itemize}
\item
 $\lambda(\Omega)>0$ if and only if $\Omega=\alpha+\beta$ where $\alpha,\beta$ are real  decomposable $3$-forms and $\alpha\wedge \beta\ne 0$
\item
$\lambda(\Omega)<0$ if and only if $\Omega=\alpha+\bar\alpha$ where $\alpha \in \Lambda^3(W^*\otimes \C)$ is a complex decomposable $3$-form and $\alpha\wedge \bar\alpha\ne 0$
\end{itemize}
\end{prp}
\begin{prf} Let $V=W\otimes \C$ be the complexification of $W$. From Proposition \ref{decomprop}, 
$\Omega=\alpha+\beta$ for decomposable complex $3$-forms. Since $\Omega$ is real,
 and $\alpha,\beta$ are unique up to ordering, complex conjugation must preserve the pair and
   there are only two possibilities: either 
$\alpha$ and $\beta$ are both real or  
$\beta=\bar\alpha.$

To decide which holds, recall that the definition of $\alpha$ and $\beta$ in the proof of the proposition was
\begin{eqnarray*}
\Omega+{\lambda(\Omega)}^{-3/2}K_{\Omega}^*\Omega&=&2\alpha\\
\Omega- {\lambda(\Omega)}^{-3/2}K_{\Omega}^*\Omega&=&2\beta
\end{eqnarray*}
Thus if $\lambda(\Omega)>0$, $\lambda(\Omega)^{1/2}$ is real and $\alpha$ and
$\beta$ are real,  and if $\lambda(\Omega)<0$, the square root  is imaginary and $\alpha$ and $\beta$ are complex conjugate. \end{prf}
We deduce from Proposition \ref{signprop} that if $\Omega$ is real and $\lambda(\Omega)>0$, then it lies in the $GL(W)$ orbit of 
$$\varphi=\theta_1\wedge \theta_2\wedge \theta_3+\theta_4\wedge \theta_5\wedge \theta_6$$
for a basis $\theta_1,\dots,\theta_6$ of $W^*$, and if $\lambda(\Omega)<0$ in the orbit of 
\begin{equation}
\varphi=\alpha+\bar\alpha,\quad\alpha=(\theta_1+i\theta_2)\wedge (\theta_3+i\theta_4)\wedge (\theta_5+i\theta_6)\label{Oo}
\end{equation}
Thus  the $20$-dimensional real vector space $\Lambda^3W^*$ contains an invariant quartic hypersurface $\lambda(\Omega)=0$ which divides $\Lambda^3W^*$ into two open sets: $\lambda(\Omega)>0$ and $\lambda(\Omega)<0$. The identity component of the stabilizer of a  $3$-form lying in the former is conjugate to $SL(3,\R)\times SL(3,\R)$ and in the latter to $SL(3,\C)$.

Choose now an orientation on $W$: this is a class of bases for  the one-dimensional vector space $\Lambda^6W^*$. We  then have a distinguished ordering of  $\alpha$ and $\beta$ or $\alpha$ and $\bar\alpha$, by the condition that $\alpha\wedge\beta$ or $i\alpha\wedge \bar \alpha$ lies in the orientation class. 
\begin{definition} Let $W$ be  oriented and  $\Omega\in \Lambda^3W^*$ be such that $\lambda(\Omega)\ne 0$. Then, writing $\Omega$ in terms of  decomposables ordered by the orientation, we define $\hat \Omega\in \Lambda^3W^*$ by
\begin{itemize}
\item
if $\lambda(\Omega)>0$, and $\Omega=\alpha+\beta$ then  $\hat\Omega=\alpha-\beta$
\item
if $\lambda(\Omega)<0$, and $\Omega=\alpha+\bar\alpha$ then  $\hat\Omega=i(\bar\alpha-\alpha)$.
\end{itemize}
\end{definition}
Note that $\hat{\hat\Omega}=-\Omega$ in both cases. The complementary $3$-form $\hat\Omega$ has the defining property that if $\lambda(\Omega)>0$ then $\Omega+\hat\Omega$ is decomposable and if $\lambda(\Omega)<0$, the complex form $\Omega+i\hat\Omega$ is decomposable.

\vskip .25cm
In this paper we shall be mainly concerned with the open set 
$$U=\{\Omega \in \Lambda^6W^*: \lambda(\Omega)<0\}$$
As we have seen, the stabilizer of $\Omega\in U$ is conjugate to $SL(3,\C)$ so $U$ is just the homogeneous space $GL^{+}(6,\R)/SL(3,\C)$. This is the space of complex structures on $\R^6$ together with a complex-linear $3$-form. Here it appears, exceptionally,  not as a homogeneous space but as an open set in a vector space.  In concrete terms, the  real $3$-form $\Omega$ determines the structure of a complex vector space with a complex $3$-form on the real vector space $W$  as follows.

From (\ref{Ksquare}) we have $K^2_{\Omega}=\lambda(\Omega)1$ so if $\lambda(\Omega)<0$, we 
define the  complex structure $I_{\Omega}$ on $W$ by
\begin{equation}
I_{\Omega}=\frac{1}{\sqrt{-\lambda(\Omega)}}K_{\Omega}
\label{complex}
\end{equation}
The real $3$-form $\Omega$ is then the real part of the complex form of type $(3,0)$
\begin{equation}
\Omega^c=\Omega+i\hat\Omega
\label{omegac}
\end{equation}
We also have
\begin{equation}
\Omega \wedge \hat\Omega=2\sqrt{-\lambda(\Omega)}
\label{func}
\end{equation}
These statements can be proved simply by checking for $\varphi$ in 
(\ref{Oo}) since $U$ is an orbit. \vskip .25cm
We have described here the linear algebra of $3$-forms in $6$ dimensions. It will be useful also to describe the same objects in {\it symplectic} terms, which we shall do next.

\section{Symplectic geometry}

\subsection{The moment map}

Choose a non-zero vector $\epsilon \in \Lambda^6W^*$, and consider 
the group of linear transformations $SL(W)\subset GL(W)$ preserving $\epsilon$.

First note that $\Lambda^3W^*$ is a symplectic vector space with the
symplectic form $\omega$ defined by 
\begin{equation} 
\omega(\Omega_1,\Omega_2)\epsilon=\Omega_1\wedge\Omega_2 \in \Lambda^6W^*
\label{symp}
\end{equation}
This is invariant under the action of $SL(W)$, which is simple, and so 
we have a well-defined moment map $\mu:\Lambda^3W^*\rightarrow \lie{sl}(W)^*$ where $\lie{sl}(W)$,
 the Lie algebra of $SL(W)$, is the space of trace zero endomorphisms of $W$. 
 If we identify the Lie algebra $\lie{sl}(W)$ with its dual using the bi-invariant 
 form $\tr(XY)$, then we can characterize the linear transformation $K_{\Omega}$ defined in (\ref{Komega}) as follows:
  \begin{prp}
\label{momentprop}
 The moment map for $SL(W)$ acting on $\Lambda^3W^*$ is given by
$$\mu(\Omega)=K_{\Omega}.$$
\end{prp}
\begin{prf} If $V$ is a symplectic vector space with symplectic form $\omega$, the Lie algebra $\lie{sp}(V)$ of $Sp(V)$ is $Sym^2(V^*)$ the space of homogeneous quadratic polynomials on $V$. This is the space of Hamiltonian functions for $\lie{sp}(V)$. Concretely, given $a\in \lie{sp}(V)$, we have the corresponding  function
\begin{equation}
\mu_a(v)=\omega(a(v),v)
\label{ham}
\end{equation}
and the moment map $\mu$ is defined by
$$\tr(\mu(v)a)=\mu_a(v)$$
In our case  $a\in \lie{sl}(W)$ defines $\rho(a)\in \lie{sp}(\Lambda^3W^*)$ via the exterior power representation   and from (\ref{symp}) and (\ref{ham})
$$\mu_a(\Omega)\epsilon=\rho(a)\Omega\wedge\Omega$$
Using a basis $w_1,\dots, w_6$ of $W$ and its dual basis $\theta_1,\dots,\theta_6$ we write 
$$a=\sum_{i,j}a^i_jw_i\otimes \theta_j$$
and then the Lie algebra action is 
\begin{equation}
\rho(a)\Omega=\sum_{i,j}a^i_j\theta_j\wedge\iota(w_i)\Omega
\label{action}
\end{equation}
so that
\begin{equation}
\rho(a)\Omega\wedge\Omega=\sum_{i,j}a^i_j\theta_j\wedge\iota(w_i)\Omega \wedge \Omega
\label{rho}
\end{equation}
But from the definition of $K_{\Omega}$ in (\ref{Komega}),
$$\theta\wedge\iota(w)\Omega\wedge \Omega=\theta(K_{\Omega}(w))\epsilon$$ for any $\theta \in W^*$ 
so from (\ref{rho}),
$$\rho(a)\Omega\wedge\Omega=\sum_{i,j}a^i_j\theta_j(K_{\Omega}(w_i))=\tr(aK_{\Omega})$$
which proves the proposition.
\end{prf}
\begin{rmk} Note that $\lambda(\Omega)$, which (by trivializing $\Lambda^6W^*$ with $\epsilon$) is now an $SL(W)$-invariant function, is just given by
 $$\lambda(\Omega)=\frac{1}{6}\tr(\mu(\Omega)^2)$$ and thus $\lambda(\Omega)$ has a natural symplectic interpretation. Its exceptional property is that $SL(W)$ acts transitively on the generic hypersurface $\tr(\mu(\Omega)^2)=const.$ If we ask this of a general symplectic representation we obtain a finite list of possibilities which appears in the classification of symplectic holonomy groups in \cite{Merk}. This is the symplectic analogue of the fact that irreducible Riemannian holonomy groups are the compact groups which act transitively on spheres. We hope  to study the analogous special geometry associated to these other cases in subsequent papers.  
\end{rmk}

\subsection{The Hamiltonian function $\phi$}

The symplectic form restricts to the open set $U$ to define the structure of a symplectic manifold. On $U$, $\lambda(\Omega)<0$, and we define
$$\phi(\Omega)=\sqrt{-\lambda(\Omega)}$$
so that $\phi$ is a smooth function here, homogeneous of degree $2$. We can then rewrite  (\ref{complex}) as 
$$K_{\Omega}=\phi I_{\Omega}$$
 and (\ref{func}) as
 \begin{equation}
 \Omega\wedge \hat \Omega=2\phi \epsilon
 \label{funcphi}
 \end{equation}
 Note that from (\ref{funcphi}) we have
 $$2\phi(\hat\Omega) \epsilon=\hat\Omega\wedge \hat{\hat \Omega}=-\hat\Omega\wedge  \Omega=\Omega\wedge \hat \Omega$$ and so
$$\phi(\hat\Omega)=\phi(\Omega).$$
The function $\phi$ defines a Hamiltonian vector field $X_{\phi}$ on $U$. Since $U$ is an open set in the vector space $\Lambda^3W^*$, we may consider the vector field as a function
$$X_{\phi}:U\rightarrow \Lambda^3W^*$$
and then we have:
\begin{prp}
\label{hamprop}
 Let $\Omega\mapsto\hat\Omega$ be the transformation of Definition 2, then $$X_{\phi}(\Omega)=-\hat\Omega.$$
\end{prp}
\begin{prf} Write $\Omega=\alpha+\bar\alpha$ in terms of complex decomposables. Then from (\ref{funcphi})
$$2\phi \epsilon=\Omega\wedge \hat \Omega=2i\alpha\wedge\bar\alpha$$
and so
\begin{equation}
\dot\phi\epsilon=i{\dot\alpha}\wedge\bar\alpha+i\alpha\wedge{\dot{\bar\alpha}}
\label{phidot}
\end{equation}
Now if $\alpha=\theta_1\wedge \theta_2\wedge \theta_3$ is a smooth curve of decomposables, on differentiation we see that
$$\dot \alpha={\dot\theta_1}\wedge \theta_2\wedge \theta_3+{\theta_1}\wedge {\dot\theta_2}\wedge \theta_3+{\theta_1}\wedge \theta_2\wedge {\dot\theta_3}$$ and so $\dot\alpha\wedge \alpha=0.$ Hence
$$\omega(i(\alpha-\bar\alpha),{\dot\alpha}+{\dot{\bar\alpha}})\epsilon=i(\alpha-\bar\alpha)\wedge({\dot\alpha}+{\dot{\bar\alpha}})=i{\dot\alpha}\wedge\bar\alpha+i\alpha\wedge{\dot{\bar\alpha}}$$
and so from (\ref{phidot}) we can write
$$\dot\phi\epsilon=\omega(i(\alpha-\bar\alpha),{\dot\alpha}+{\dot{\bar\alpha}})\epsilon
=-\omega(\hat\Omega,\dot\Omega)\epsilon$$

Since the Hamiltonian vector field is defined by the property $\iota(X_{\phi})\omega=d\phi$ we have our result.
\end{prf}
\begin{rmk}
If we regard $U$ as parametrizing complex structures on $W$ together with
$(3,0)$ forms then 
we have an obvious circle action $\Omega+i\hat\Omega \mapsto
e^{i\theta}(\Omega+i\hat\Omega)$. The real point of view we are following here 
describes this as an action on the real parts which is 
 $$\Omega \mapsto \cos{\theta}\, \Omega-\sin{\theta}\,\hat \Omega$$
  and differentiating at $\theta=0$ we see that this action is generated 
  by the vector field $X_{\phi}=-\hat\Omega$. The position vector (or Euler vector field) 
  $\Omega$ integrates to give the $\R^*$-action $\Omega \mapsto \lambda\Omega$ and the two define the action of multiplying the complex form $\Omega+i\hat\Omega$ by $\C^*$. \end{rmk}

\subsection{The complex structure on $U$}

\begin{prp}
\label{integprop}
 The derivative $J=DX_{\phi}\in \End(\Lambda^3W^*)$ of $X_{\phi}:U\rightarrow \Lambda^3W^*$ at $\Omega$ satisfies $J^2=-1$ and defines an integrable complex structure on $U$.
\end{prp}
\begin{prf} We have seen that $\hat{\hat\Omega}=-\Omega$ and $X_{\phi}=-\hat\Omega$, so it follows that $(DX_{\phi})^2=-1$. This is an almost complex structure.

To prove integrability, it is easiest to use flat coordinates $x_1,\dots,x_{20}$ on $U$. 
The symplectic form can be written 
$\omega=\sum_{\alpha,\beta}\omega_{\alpha \beta}dx_{\alpha}\wedge dx_{\beta}$ 
with $\omega_{\alpha\beta}$ constant. Let $\omega^{\alpha\beta}$ be the inverse matrix. 
In these coordinates, since $X_{\phi}$ is the Hamiltonian vector field of
$\phi$, 
the matrix of $J=DX_{\phi}$ is \begin{equation}
J^{\alpha}_{\beta}=\sum_{\gamma}\omega^{\alpha\gamma}\frac{\partial^2 \phi}{\partial x_{\gamma}\partial x_{\beta}}
\label{Jmatrix}
\end{equation}
Following \cite{Hit}, we define the complex functions $z_1,\dots,z_{20}$ by
$$z_{\alpha}=x_{\alpha}-i\sum_{\beta} \omega^{\alpha\beta}\frac{\partial \phi}{\partial x_\beta}$$
Then
$$dz_{\alpha}=dx_{\alpha}-i\sum_{\beta,\gamma} \omega^{\alpha\beta}\frac{\partial^2 \phi}{\partial x_{\beta}\partial x_{\gamma}}dx_{\gamma}=dx_{\alpha}-i\sum_lJ^{\alpha}_{\gamma}dx_{\gamma}$$
Now
$$J(dz_{\alpha})=\sum_{\beta} J^{\alpha}_{\beta}(dx_{\beta}-i\sum_lJ^{\beta}_{\gamma}dx_{\gamma})=\sum_{\beta} J^{\alpha}_{\beta}dx_{\beta}+idx_{\alpha}=idz_{\alpha}$$
so $dz_{\alpha}$ is of type $(1,0)$. Since $dz_{\alpha}+d\bar z_{\alpha}=2dx_{\alpha}$ spans the cotangent space at each point, we can find $10$  independent functions amongst these which are local holomorphic coordinates for  $U$.
\end{prf}
 At each $\Omega$ we have the natural complex structure $J$ on the vector space $\Lambda^3W^*$ but $W$ itself is also complex with respect to $I_{\Omega}$. Take the type decomposition of $\Lambda^3W^*\otimes \C$ with respect to $I_{\Omega}$:
$$\Lambda^3W^*\otimes \C=\Lambda^{3,0}\oplus\Lambda^{2,1}\oplus\Lambda^{1,2}\oplus\Lambda^{0,3}.$$
 \begin{prp} 
 \label{Jprop}
 The complex structure $J=i$ on $\Lambda^{3,0}\oplus\Lambda^{2,1}$ and $-i$ on $\Lambda^{1,2}\oplus\Lambda^{0,3}$.
 \end{prp}
 \begin{prf} Given $a\in \lie{sl}(W)$ let $X_a$ be the vector field on $\Lambda^3W^*$ induced by the action. This is a linear action. This means that   $DX_a=\rho(a)$ where  $\rho$ is the representation at the Lie algebra level, and $X_a=\rho(a)\Omega$ where $\Omega$ is the position vector.

 Now $\phi$ is invariant under $SL(W)$ so $[X_a,X_\phi]=0$, thus
\begin{equation}
J(X_a)=DX_{\phi}(X_a)=DX_a(X_{\phi})=\rho(a)X_{\phi}
\label{JJ}
\end{equation}
But $X_{\phi}=-\hat\Omega$ from Proposition \ref{hamprop}, and  $X_a=\rho(a)\Omega$, thus (\ref{JJ}) can be written
$$J(\rho(a)\Omega)=-\rho(a)\hat\Omega$$
and since $J^2=-1$,
\begin{equation}
J(\rho(a)(\Omega+i\hat\Omega))=-\rho(a)(\hat\Omega-i\Omega)=i\rho(a)(\Omega+i\hat\Omega)
\label{j2}
\end{equation}
From (\ref{action})
$$\rho(a)(\Omega+i\hat\Omega)=\sum_{i,j}a^i_j\theta_j\wedge\iota(w_i)(\Omega+i\hat\Omega)$$
and since $\Omega+i\hat\Omega$ is of type $(3,0)$, this is of type $(3,0)+(2,1)$. 
The complex structure $J$ acts on it as $i$ from (\ref{j2}). 

Consider $J$ acting on the  vector field $\Omega$. We  use (\ref{Jmatrix}) to calculate
$$J(\Omega)_{\alpha}=\sum_{\beta,\gamma}\omega^{\alpha\gamma}x_{\beta}\frac{\partial^2 \phi}{\partial x_{\gamma}\partial x_{\beta}}=\sum_{\gamma}\omega^{\alpha\gamma}\frac{\partial \phi}{\partial x_{\gamma}}=-\hat\Omega_{\alpha}$$
since $\phi(\Omega)$ is homogeneous of degree $2$ and so $\partial
\phi/\partial x_{\alpha}$ has degree 1. It follows that
$J(\Omega+i\hat\Omega)=i(\Omega+i\hat\Omega)$ so $J$ acts as $i$ on forms of type $(3,0)$. 

Now $SL(W)$ has orbits 
of real codimension one, and the position vector field $\Omega$ is transverse to the orbit so any vector of type
$(3,0)+(2,1)$ is a complex linear combination of a vector  $\rho(a)(\Omega+i\hat\Omega)$ and $\Omega+
i\hat\Omega$ and hence is acted on by $J$ as $i$. \end{prf}

\subsection{Special geometry of $U$}

 We have actually derived here a certain geometric structure  on the open set $U$ --  a {\it special pseudo-K\"ahler metric}. Recall from \cite{Freed} the definition: \begin{definition} A special K\"ahler manifold is a complex manifold $M$ with complex structure $J\in \Omega^1(T)$ such that
\begin{itemize}
\item
there is a K\"ahler metric $g$ with K\"ahler form $\omega$
\item
a flat torsion-free connection $D$ such that 
\item
$D\omega=0$ and
\item
$d_{D}J=0\in \Omega^2(M,T)$
\end{itemize}
\end{definition}
\begin{prp} 
\label{openprop}
The open set $U\subset \Lambda^3W^*$ has an $SL(W)$-invariant special pseudo-K\"ahler structure of Hermitian signature $(1,9)$ and $X_{\phi}$ is an infinitesimal isometry of the metric.
\end{prp}
\begin{prf}  Firstly, $U$ is given as an open set in a symplectic vector space, so we take the flat torsion-free connection as the ordinary derivative $D$. The complex structure $J$ is given by $J=DX_{\phi}=d_{D}X_{\phi}$, so 
$$d_{D}J=d_{D}^2X_{\phi}=0$$
since $D$ is flat. For a pseudo-K\"ahler metric, the metric and symplectic form are related by
$$g(X,Y)=\omega(JX,Y).$$
But  $\omega(J\Omega_1,\Omega_2)$ is given in flat coordinates by
$$\sum_{\alpha} J^{\alpha}_{\beta}\omega_{\alpha\gamma}=\frac{\partial^2 \phi}{\partial x_{\gamma}\partial x_{\beta}}$$
from (\ref{Jmatrix}). This is symmetric and hence defines a metric
$$g=\sum_{\beta,\gamma} \frac{\partial^2 \phi}{\partial x_{\beta}\partial x_{\gamma}}dx_{\beta}dx_{\gamma}.$$
 We need  to determine the signature of the corresponding Hermitian form.

From Proposition \ref{Jprop}, the $(1,0)$ vectors for $J$ are $\Lambda^{3,0}\oplus\Lambda^{2,1}$ and the symplectic form $\omega$ is 
$$\omega(\Omega_1,\Omega_2)\epsilon=\Omega_1\wedge\Omega_2.$$
Now if $\theta_1,\theta_2,\theta_3$ form a basis for $\Lambda^{1,0}$,
$$\theta_1\wedge \theta_2\wedge \bar\theta_3\wedge \bar\theta_1\wedge \bar 
\theta_2\wedge \theta_3=-\theta_1\wedge \theta_2\wedge \theta_3\wedge \bar
\theta_1\wedge \bar\theta_2\wedge \bar\theta_3$$ so the Hermitian form on
$\Lambda^{2,1}$ has the opposite sign to that on the $1$-dimensional space $\Lambda^{3,0}$. But the 
metric applied to the position vector field $\Omega$ (which is the real part of a
vector in $\Lambda^{3,0}$) is  $$g(\Omega,\Omega)=\sum\frac{\partial^2
\phi}{\partial x_{\alpha}\partial x_{\beta}}x_{\alpha}x_{\beta}=2\phi$$ since $\phi$ is homogeneous of degree $2$. Since $\phi$ is positive, the metric is positive on $\Lambda^{3,0}$ and hence has signature $(1,9)$.

The vector field $X_{\phi}$ integrates, as we have seen, to the circle action $\Omega+i\hat\Omega \mapsto e^{i\theta}(\Omega+i\hat\Omega)$. This changes the complex $3$-form, but not the underlying complex structure $I_{\Omega}$. From Proposition \ref{Jprop}, the complex structure $J$ on $U$ only depends on $I_{\Omega}$, so that the circle action preserves $J$. By  definition $X_{\phi}$ is symplectic, so it preserves the metric too.
\end{prf}

\begin{rmk} The $\C^*$ action generated by the vector fields
$X_{\phi}=\hat\Omega$ and  $\Omega$ is holomorphic and
the quotient has (see \cite{Freed}) a {\it projective special K\"ahler structure}. This space in our case is the space of complex structures on $\R^6$ compatible with a given orientation, the homogeneous space $GL^{+}(6,\R)/GL(3,\C)$. \end{rmk}
 The symplectic approach of this section will be especially useful in analyzing the variational problem, but we take a short detour first to investigate a little more the linear algebra of $3$-forms in $6$ dimensions.

\section{Self-duality of $3$-forms}

\subsection{Lorentzian self-duality}

We  point out here how the  geometry of $3$-forms in six dimensions can be interpreted in the
presence of a Lorentzian metric. This will shed some light on the nonlinear map
$\Omega\mapsto \hat\Omega$ and also explain the setting for a nonlinear
equation of some current interest in theoretical physics.

Suppose then that the $6$-dimensional real vector space $W$ has a metric of signature $(5,1)$. The Hodge star operator  $$\ast: \Lambda^3 W^*\rightarrow \Lambda^3 W^*$$ is defined by the basic property
\begin{equation}
\alpha\wedge\ast \beta=(\alpha,\beta)\epsilon
\label{hodge}
\end{equation}
where $\epsilon \in \Lambda^6 W^*$ is the volume form. 
\begin{ex}
If $e_0,e_1,\dots,e_5$ is an orthogonal basis of $W^*$ with $(e_0,e_0)=-1$ and  $(e_i,e_i)=1$ for $1\le i\le 5$, then with $\epsilon=e_0\wedge e_1\wedge \dots \wedge e_5$, we have 
\begin{eqnarray*}
\ast e_0\wedge e_1\wedge e_2&=&-e_3\wedge e_4\wedge e_5\\
\ast e_3\wedge e_4\wedge e_5&=&-e_0\wedge e_1\wedge e_2
\end{eqnarray*}
\end{ex}
The star operator just defined has the property $\ast^2=1$ and we can decompose $\Lambda^3W^*$ into the $\pm 1$ eigenspaces of $\ast$ -- the self-dual and anti-self-dual $3$-forms:
$$\Lambda^3W^*=\Lambda_+\oplus \Lambda_-$$
If $\beta =\theta_1\wedge\theta_2\wedge \theta_3$ is a non-zero decomposable $3$-form, 
let $\varphi_1,\varphi_2,\varphi_3$ be  a basis for the $3$-dimensional space orthogonal
 to that spanned by $\theta_1,\theta_2, \theta_3$. Then from the definition of the Hodge star
  (\ref{hodge}), we see that $\alpha\wedge\ast \beta=0$ for any $3$-form $\alpha$ of 
  the form $\varphi_i\wedge \rho$. This means that $\ast \beta$ is a multiple of 
  $\varphi_1\wedge\varphi_2\wedge \varphi_3$. Thus, as in the example above, 
   the Hodge star takes {\it decomposable $3$-forms to decomposable $3$-forms}.
   Geometrically it transforms a volume form on a $3$-dimensional subspace to  a volume form on the orthogonal space.

If $\beta$ (real or complex) is decomposable and $\ast \beta=\beta$ then the space spanned by $\theta_1,\theta_2, \theta_3$ is orthogonal to itself and hence isotropic. 

We see next that real self-dual $3$-forms lie on one side of the hypersurface $\lambda(\Omega)=0$:
\begin{prp} 
\label{sdprop}
Let $\Omega \in \Lambda^3 W^*$ be self-dual then $\lambda(\Omega)\ge 0$ and if $\lambda(\Omega)> 0$ then $\Omega=\alpha+\ast\alpha$ where $\alpha \in \Lambda^3 W^*$ is decomposable.
\end{prp}
\begin{prf} Suppose $\lambda(\Omega)< 0$, then from Proposition \ref{decomprop},  $\Omega=\alpha +\bar\alpha$ with $\alpha=\theta_1\wedge\theta_2\wedge \theta_3$ and $\theta_i$ complex. These span the $(1,0)$ forms for the complex structure $I_{\Omega}$.  Since $\ast$ is linear and $\alpha,\bar\alpha$ are unique up to ordering, the self-duality condition implies  that either
$${\ast\alpha=\alpha}\quad{\rm or}\quad {\ast\alpha=\bar\alpha}.$$
In the first case, the space of $(1,0)$ forms is isotropic, so we have a 
 pseudo-hermitian metric. This must have real signature of the form $(2k,6-2k)$
 and  not $(5,1)$.

In the second case the $\theta_i$ are orthogonal to the $\bar\theta_j$ so the metric is of type $(2,0)+(0,2)$ and so the signature must be $(3,3)$ and not $(5,1)$.

We conclude that $\lambda(\Omega)>0$ and from Proposition \ref{decomprop},
$\Omega=\alpha+\beta$ for real decomposable forms $\alpha,\beta$. Since $\ast \Omega=\Omega$ we must have either $\ast \alpha=\alpha$ and $\ast \beta=\beta$ or $\ast \alpha=\beta$. The first implies that there is a real isotropic $3$-dimensional subspace which is impossible in signature $(5,1)$ so the second condition holds, as required. 

Note that since $\alpha\wedge \beta\ne 0$, the two $3$-dimensional subspaces of $W$ determined by $\alpha$ and $\beta$ span $W$. Since they are orthogonal, one must have signature $(2,1)$ and the other $(3,0)$.
\end{prf}

\begin{rmk}
We see from the proposition that $\Omega=\alpha+\ast\alpha$ where $\alpha$ and $\ast\alpha$ are decomposable. The signature of the two $3$-dimensional subspaces gives us a means of distinguishing $\alpha$ and $\ast \alpha$ and hence an ordering. The formula (\ref{hodge}) tells us that
$$\alpha\wedge \ast \alpha=(\alpha,\alpha)\epsilon$$
so if the metric on the $3$-dimensional subspace of $W$ determined by $\alpha$ is of signature $(2,1)$ (or equivalently $(\alpha,\alpha)>0$) then this is the ordering determined as in Section 2 by the orientation.
\end{rmk}

Since for a self-dual form $\lambda(\Omega)>0$, and $\Omega=\alpha+\ast\alpha$, by 
the definition of the map $\Omega\mapsto \hat\Omega$ we have 
$$\hat\Omega=\alpha-\ast\alpha$$
hence
\begin{prp}
\label{sdual}
If $\Omega$ is self-dual $\hat\Omega$ is  anti-self-dual.
\end{prp}
 Thus $\Omega\mapsto \hat\Omega$ is 
a nonlinear map from an open set of the space $\Lambda_+$ to $\Lambda_-$. We define
 in this case $$\phi(\Omega)=\sqrt{\lambda(\Omega)}$$ to be the positive square
 root. 
 
 \subsection{Lagrangian aspects}
 
There is a symplectic interpretation also to  self-duality. 
First recall that the symplectic form $\omega$ on $\Lambda^3W^*$ is defined by
$\omega(\Omega_1,\Omega_2)\epsilon=\Omega_1\wedge\Omega_2$. We have
\begin{prp}
\label{lagprop} The subspaces $\Lambda_+, \Lambda_-\subset \Lambda^3W^*$ of self-dual and anti-self-dual $3$-forms are Lagrangian.
\end{prp}
\begin{prf}
Suppose  $\Omega_1,\Omega_2$ are self-dual then
$$\Omega_1\wedge\Omega_2=\Omega_1\wedge \ast\Omega_2=(\Omega_1,\Omega_2)\epsilon$$
 and similarly $\Omega_2\wedge\Omega_1=(\Omega_2,\Omega_1)\epsilon$. But $(\Omega_1,\Omega_2)=(\Omega_2,\Omega_1)$ and 
$\Omega_1\wedge\Omega_2=-\Omega_1\wedge\Omega_2$ so 
$\Omega_1\wedge\Omega_2=0.$
The same argument gives the anti-self-dual version.
\end{prf}
From this result, the symplectic pairing identifies $\Lambda_-$ with the dual space $\Lambda_+^*$.  We then obtain:
\begin{prp}
\label{hatprop}
The map $\Omega \mapsto \hat\Omega$ from $\Lambda_+$ to $\Lambda_-\cong 
\Lambda_+^*$ is the derivative of $-\phi$ restricted to $\Lambda_+$. \end{prp}
\begin{prf} We saw in Proposition \ref{hamprop} that $-\hat \Omega=X_{\phi}$ the Hamiltonian vector field generated by the function $\phi$. Let $\dot\Omega$ be tangent to $\Lambda_+$, then
$$d\phi(\dot\Omega)=\omega(X_{\phi},\dot\Omega)$$
But this pairing identifies $\Lambda_-$, in which $X_{\phi}=-\hat\Omega$ lies, with $\Lambda_+^*$, hence the result.
\end{prf}

From (\ref{lagprop}) and (\ref{hatprop}) we see that the subset 
$$L=\{\Omega +\hat\Omega\in \Lambda_+\oplus \Lambda_-:  \lambda(\Omega)> 0\}$$
 is the  Lagrangian submanifold generated by the invariant function $-\phi$. If
 $f(x)$ is an arbitrary smooth function, then the form $$\Omega+f(\phi)\hat\Omega$$
where $\Omega$ is self-dual, lies on the Lagrangian submanifold generated by $g(\phi)$
where $g'(x)=f(x)$. This is the setting  for a nonlinear equation  -- the equation of motion for a self-interacting self-dual tensor in six dimensions -- which appears in the physics literature \cite{HSW} in the context of M-theory five-branes. Explicitly these equations are:
$$dH=d(H_{+} + H_{-})=0$$
where $H_{+},H_{-}$ are the self-dual and anti-self-dual components of the $3$-form $H$, and
\begin{eqnarray*}
(H_{+})_{abc}&=&Q^{-1}h_{abc}\\
(H_{-})_{abc}&=&Q^{-1}k_a^dh_{dbc}
\end{eqnarray*}
where
$$k_a^b=h_{acd}h^{bcd}$$
and 
$$Q=1-\frac{2}{3}\tr k^2.$$
The equations say in our language that $H_{-}$ is proportional to $\hat H_{+}$. 
They  can be put in the form
$$d(\Omega+f(\phi)\hat\Omega)=0$$
for a particular function $f(\phi)$. 

If we replace the nonlinear Lagrangian submanifold which defines 
this equation by the linear one $\Lambda_+$ then we just obtain the self-dual ``Maxwell equations"  in six dimensions.

\subsection{Spinor formulation}

For self-dual forms, the function $\phi$ and the nonlinear map $\Omega \mapsto \hat\Omega$ 
also have a concrete representation using $2$-component spinors. We recall the special isomorphism 
$$Spin(5,1)\cong SL(2,\K).$$ 
The two spin representations in this signature are dual $4$-dimensional complex vector spaces $S$ and $S^*$ with a quaternionic structure, and then
$$\Lambda_+\cong Sym^2 S^*,\quad \Lambda_-\cong Sym^2 S$$
are the real $10$-dimensional spaces of self-dual and anti-self-dual $3$-forms. Thus a self-dual $3$-form can be interpreted as a symmetric linear map
$$A:S\rightarrow S^*.$$
The quartic invariant function  $\lambda$ must be a multiple of the only $SL(4,\C)$  quartic invariant of $A$,
 namely the determinant:
$$\lambda(A)=4 \det A.$$
We find now the following simple expression for $\hat\Omega$:
\begin{prp}
\label{dualprop}
\begin{itemize}
\item
 If $A\in Sym^2 S^*$ represents a self-dual $3$-form $\Omega$ then 
 $$\hat\Omega=\sqrt{\det A}A^{-1}\in Sym^2 S.$$
 \item
 If $B\in Sym^2 S$ represents an anti-self-dual $3$-form $\Omega$ then 
 $$\hat\Omega=-\sqrt{\det A}A^{-1}\in Sym^2 S.$$
 \end{itemize}
\end{prp}
\begin{prf}
We use Proposition \ref{hatprop}. We have the standard identity for differentiating the determinant: 
$$\frac{1}{\det A}d(\det A)(\dot A)=d(\log \det A)(\dot A)=\tr(A^{-1}\dot A)$$
and so since $\phi=\sqrt{\lambda}=2\sqrt{\det A}$,
$$d\phi(\dot A)=d(2\sqrt{\det A})(\dot A)=\tr (\sqrt{\det A}A^{-1}\dot A)$$
and taking $\tr(AB)$ to be the natural pairing between $Sym^2 S$ and $Sym^2 S^*$ we have the result for self-dual forms. The second result follows from $\hat{\hat\Omega}=-\Omega$.
\end{prf}

 The isomorphism $Spin(5,1)\cong SL(2,\K)$ is fundamental in the twistor theory of the $4$-sphere. 
 The group $SO(5,1)$ acts as the conformal transformations of $S^4$ defined by the quadric 
 $x_0^2=x_1^2+\dots+x_5^2$ in $\RP^5=P(W)$ and $SL(2,\K)$ as the projective transformations of the 
 twistor space $\CP^3=P(S)$ which commute with the real structure. The $4$-sphere parametrizes this way 
 the real lines in $\CP^3$. Now a self-dual $3$-form  with $\lambda(\Omega)>0$ can be
  written $\Omega =\alpha +\ast\alpha$ where $\alpha$ defines a $3$-dimensional subspace of 
  $W$ with signature $(2,1)$. The corresponding plane in $\RP^5$  intersects $S^4$ in a  circle. 
  On the other hand the quadratic form $A$ representing $\Omega$ defines a quadric in $\CP^3$. 
  This geometry comes to the fore in particular in the study of the moduli space of charge $2$ instantons \cite{Hart}.

Note that in the twistor interpretation, a quadratic form $A$ with $\det(A)\ne 0$ defines a nonsingular quadric in $\CP^3$. From (\ref{dualprop}) the map $\Omega\mapsto \hat \Omega$ is then nothing more than replacing  a quadric by the dual quadric.

\section{An invariant functional}

\subsection{Critical points}

We return now to the mainstream of our development. Let $M$ be a closed, oriented $6$-manifold. 
A global $3$-form $\Omega$ on $M$ defines at each point a vector in $\Lambda^3T^*$ and so  
a global section $\lambda(\Omega)$ of $(\Lambda^6T^*)^2$. We take $\vert
\lambda(\Omega)\vert$, which is a 
non-negative continuous  section of $(\Lambda^6T^*)^2$, and the  square root (a
section of $\Lambda^6 T^*$) which lies in the given orientation class. Thus $\sqrt{\vert \lambda(\Omega)\vert}$ is a continuous $6$-form on $M$, which is smooth wherever $\lambda(\Omega)$ is non-zero. We define a functional $\Phi$ on $C^{\infty}(\Lambda^3T^*)$ by  \begin{equation}
\Phi(\Omega)=\int_M \sqrt{\vert \lambda(\Omega)\vert}
\label{functional}
\end{equation}
Since $\lambda:\Lambda^3T^*\rightarrow  (\Lambda^6T^*)^2$ is $GL(6,\R)$-invariant,  $\Phi$ is 
clearly invariant under the action of orientation-preserving diffeomorphisms. 
\vskip .25cm
The functional $\Phi$ is homogeneous of degree $2$ in $\Omega$, just like the norm square 
of a form using  a Riemannian metric and it is natural to try and apply
``non-linear Hodge theory" -- to look for critical points of $\Phi$ on a
cohomology class of closed $3$-forms on $M$. We find  the following:
 \begin{thm} 
\label{CYthm}
Let $M$ be a compact complex $3$-manifold with trivial canonical bundle and
$\Omega$ the real part of a non-vanishing holomorphic $3$-form. Then $\Omega$ is a critical point of the functional $\Phi$ restricted to the cohomology class $[\Omega]\in H^3(M,\R)$.

 Conversely, if $\Omega$ is a critical point on a cohomology class of an
 oriented closed $6$-manifold $M$ 
 and $\lambda(\Omega)<0$ everywhere, then $\Omega$ defines  on $M$ the structure
 of  a  complex manifold  such that $\Omega$ is the real part of a non-vanishing holomorphic $3$-form. \end{thm}
\begin{prf} Choose a non-vanishing $6$-form $\epsilon$ on $M$ and assume that $\Omega$ is a closed $3$-form with  $\lambda(\Omega)<0$ (as is the case for the real part of a holomorphic $3$-form). Writing 
$$\lambda(\Omega)=-\phi(\Omega)^2\epsilon^2$$
we have
$$\Phi(\Omega)=\int_M\phi(\Omega)\epsilon.$$
Since $\lambda(\Omega)\ne 0$, $\phi$ is smooth so  take the first variation of this functional:
\begin{equation}
\delta \Phi(\dot \Omega)=\int_Md\phi(\dot\Omega)\epsilon
\label{1var}
\end{equation}
 From the symplectic interpretation, 
$d\phi(\dot\Omega)=\omega(X_{\phi},\dot\Omega)$
and $X_{\phi}=-\hat\Omega$, so that
$$d\phi(\dot\Omega)\epsilon=-\hat \Omega \wedge \dot \Omega$$
and therefore (\ref{1var}) gives
\begin{equation}
\delta \Phi(\dot \Omega)=-\int_M\hat \Omega \wedge \dot \Omega
\label{firstv}
\end{equation}
We are varying in a fixed cohomology class, so $\dot \Omega=d\varphi$ for some
$2$-form $\varphi$. Putting this in the integral, we have 
$$\delta \Phi(\dot \Omega)=-\int_M\hat \Omega \wedge d\varphi=\int_M d\hat\Omega \wedge \varphi$$
by Stokes' theorem. Thus $\delta\Phi=0$ at $\Omega$ for all $\varphi$ if and only if
$$d\hat\Omega=0.$$ 
If $M$ is a complex manifold with a non-vanishing holomorphic $3$-form $
\Omega+i\hat\Omega$  then since a holomorphic $3$-form is closed,  
$d\Omega=d\hat\Omega=0$ and so $\Omega$ is a critical point for $\Phi$. 

Conversely, assume that $d(\Omega+i\hat\Omega)=0$ and $\lambda(\Omega)<0$. Then  a complex $1$-form $\theta$ is of type $(1,0)$ with respect to the almost complex structure $I_{\Omega}$ if and only if 
$$\theta\wedge (\Omega+i\hat\Omega)=0.$$
Taking the exterior derivative
$$d\theta\wedge (\Omega+i\hat\Omega)=0$$
which means that $d\theta$ has no $(0,2)$ component. But from the Newlander-Nirenberg theorem this means that $I_{\Omega}$ is integrable. Since $(\Omega+i\hat\Omega)$ is of type $(3,0)$,
$$0=d(\Omega+i\hat\Omega)=\bar\partial (\Omega+i\hat\Omega)$$
and so is holomorphic.
\end{prf}
\begin{ex} A Calabi-Yau threefold is by definition a K\"ahler manifold with a covariant constant (and hence nonvanishing) holomorphic $3$-form, and so its complex structure appears as a critical point of the functional. These are perhaps the most interesting and plentiful examples. On the other hand, at the opposite extreme from the K\"ahler case are the non-K\"ahler complex threefolds with trivial canonical  bundle which are diffeomorphic to connected sums of $k\ge 2$ copies of $S^3\times S^3$ (see \cite{Tian}).
\end{ex}

\subsection{Nondegeneracy}

Since the functional $\Phi$ is diffeomorphism invariant, any critical point lies on a
$\Diff(M)$-orbit of critical points and so cannot be  non-degenerate. We can
ask however  if it is formally a Morse-Bott critical point i.e. if its Hessian is nondegenerate transverse to the action of $\Diff(M)$. The next proposition gives conditions under which this is true: \begin{prp} 
\label{nondegenprop}
Let $M$ be a compact complex $3$-manifold with non-vanishing holomorphic $3$-form $\Omega+i\hat\Omega$. Suppose $M$ satisfies the $\partial \bar\partial$-lemma, then the Hessian of $\Phi$ is nondegenerate transverse to the action of $\Diff(M)$ at $\Omega$.
\end{prp}
(Recall that the  $\partial \bar\partial$-lemma holds if each exact $p$-form
$\alpha$ satisfying $\partial\alpha=0$ and $\bar\partial\alpha=0$ can be written as 
$\alpha= \partial \bar\partial \beta$ for some $(p-2)$-form $\beta$. This is true for any K\"ahler
 manifold, but also for the non-K\"ahler examples above.)

\begin{prf} 
We need the second variation of the functional 
$$\Phi(\Omega)=\int_M\phi(\Omega)\epsilon.$$
Since $\phi$ depends only on the value of $\Omega$ at a point and not its derivatives we obtain the Hessian
$$\delta^2\Phi(\dot\Omega_1,\dot\Omega_2)=\int_M D^2\phi(\dot\Omega_1,\dot\Omega_2)\epsilon$$
but from the Hamiltonian interpretation we can write this as
\begin{equation}
\delta^2\Phi(\dot\Omega_1,\dot\Omega_2)=\int_M DX_\phi(\dot\Omega_1)\wedge \dot\Omega_2=\int_M J\dot\Omega_1\wedge \dot\Omega_2
\label{var}
\end{equation}
We want to study the degeneracy of the Hessian, so suppose that $\delta^2\Phi(\dot\Omega_1,\dot\Omega_2)=0$ for $\dot\Omega_1=d\psi$ and {\it all} $\dot\Omega_2=d\varphi$. Then from (\ref{var})
$$0=\int_M Jd\psi\wedge d\varphi=-\int_M dJd\psi\wedge \varphi$$
by Stokes' theorem. This holds for all $2$-forms $\varphi$, which implies
$$dJd\psi=0.$$
Thus to prove the theorem we should deduce from this that $d\psi$ is tangent to  the  $\Diff(M)$ orbit through $\Omega$.
\vskip .25cm 
The tangent space to the $\Diff(M)$ orbit  is the space spanned by ${\cal L}_X\Omega$ for a vector field $X$. We have
$${\cal L}_X\Omega=d(\iota(X)\Omega)+\iota(X)d\Omega$$
but since  $\Omega$ is closed, these are the $3$-forms $d(\iota(X)\Omega)$. Since $\Omega$ is the real part of a non-vanishing $(3,0)$ form, the $2$-forms $\iota(X)\Omega$ are precisely the real sections $\alpha$ of $\Lambda^{2,0}\oplus\Lambda^{0,2}$.

The theorem follows from the following lemma:
\begin{lem}
\label{ddlem}
 Let $M$ be a compact complex $3$-manifold with non-vanishing holomorphic $3$-form which satisfies the $\partial \bar\partial$-lemma, then
$d\psi$ is tangent to  the $\Diff(M)$ orbit of $\Omega$ if and only if $dJd\psi=0$.
\end{lem}
\begin{lemprf} If $\psi$ is of type $(2,0)$ then $d\psi \in C^{\infty}(\Lambda^{3,0}\oplus \Lambda^{2,1})$ and then from Proposition \ref{Jprop}, $dJd\psi=d(id\psi)=id^2\psi=0$. Thus $dJd\psi=0$.

Conversely, suppose $dJd\psi=0$. Since we have just seen this holds for all $\psi$ of type $(2,0)+(0,2)$ assume $\psi$ has type $(1,1)$.
Then  $d\psi=\partial \psi+\bar\partial \psi$ where $\partial \psi, \bar\partial \psi$ are of type $(2,1),(1,2)$ respectively. From Proposition \ref{Jprop},
$$0=dJd\psi=id\partial \psi-id\bar\partial \psi=2i\partial \bar\partial \psi.$$
This means that $\partial d\psi=\bar \partial d\psi=0$ and so applying the $\partial \bar\partial$-lemma to $d\psi$ we can write 
$$d\psi=i\partial \bar\partial\gamma$$ for a real $1$-form $\gamma$. Writing $\gamma=\theta +\bar \theta$ where $\theta$ is of type $(1,0)$, we have

$$d\psi=i\partial \bar\partial\gamma=i\partial \bar\partial(\theta + \bar \theta)=d(i(\bar\partial\bar\theta-\partial\theta))$$

Since $i(\bar\partial\bar\theta-\partial\theta)$ is real and of type $(2,0)+(0,2)$ this proves the lemma.
\end{lemprf}
\end{prf}
At a formal level what we have proved here is that if we consider 
 the invariant functional $\Phi$ as a function on the quotient of a cohomology class by $\Diff(M)$, 
 then it has a non-degenerate  critical point at a $3$-form $\Omega$ which defines 
 a complex manifold with trivial canonical bundle. We  then expect that nearby 
 cohomology classes will also have non-degenerate critical points and that an open set
  in $H^3(M,\R)$ will parametrize the moduli of such complex structures, (together with holomorphic $3$-forms). 
  This is true in all dimensions by the results of A.~Todorov \cite{Tod} and G.~Tian  \cite{Tian1}, 
  but we shall give next a direct treatment in three dimensions from our variational point of view.

 \section{The moduli space}
 
 \subsection{Sobolev spaces}
 
 If $f(x,t)$ is a smooth  family of functions $f:\R^m\times \R^n \mapsto \R$ such that $f(x,0)$ has a non-degenerate critical point at $x=0$, then there is a neighbourhood $U\times V$ of $(0,0)$ such  for each $t\in V$, $f(x,t)$ has a unique non-degenerate critical point in $U$. To prove this, we just apply the implicit function theorem to the map
$$(x,t)\mapsto D_xf.$$
 We shall use this argument next in a Banach space context to translate the formal results of the last section into a concrete construction of  a moduli space.

Take a $3$-form $\Omega$ on $M$ defining a complex structure satisfying the $\partial \bar\partial$-lemma, and choose a Hermitian metric. We take a slice for the $\Diff(M)$ action at $\Omega$ by looking at the space of closed forms which are orthogonal to the orbit of $\Omega$. Recall that the tangent space of the orbit consists of forms $d\psi$ where $\psi$ is real and of type $(2,0)+(0,2)$, so   using 
$$\int_M(\alpha,d\psi)=\int_M(d^*\alpha,\psi)$$
 orthogonality is the condition 
$$(d^*\alpha)^{2,0}=0.$$
We shall work with Sobolev spaces of forms $L^2_k(\Lambda^p)$, choosing $k$ appropriately when required. So let $E$ be the Banach space
$$E=\{\alpha \in L^2_k(\Lambda^3): d\alpha=0\quad{\rm and}\quad  
(d^*\alpha)^{2,0}=0\}.$$ First we show that  $L^2$ orthogonal  projection onto $E$ is well-behaved in the Sobolev norm.

Let $G$ be the Green's operator for the Laplacian $\Delta$ on forms. Then elliptic regularity says that 
$$G:L^2_k(\Lambda^p)\rightarrow L^2_{k+2}(\Lambda^p)$$
and for any form $\alpha$ we have a Hodge decomposition
\begin{equation}
\alpha=H(\alpha)+d(d^*G\alpha)+d^*(Gd\alpha)
\label{Hodge}
\end{equation}
where $H(\alpha)$ is harmonic. So given $\alpha \in L^2_k(\Lambda^3)$ define first the  form $\alpha_1\in L^2_k(\Lambda^3)$ by
$$\alpha_1=H(\alpha)+d(d^*G\alpha).$$
This is an $L^2$ orthogonal projection to closed forms in the same Sobolev
space. To further project onto $E$  we  want to find a form $\theta$ of type $(2,0)$ such that
$$\alpha_2=\alpha_1+d(\theta +\bar\theta)$$
satisfies $(d^*\alpha_2)^{2,0}=0$. Write 
$$d^*G\alpha=\rho+\nu+\bar\nu \in L^2_{k+1}(\Lambda^2)$$
where $\rho$ is of type $(1,1)$ and $\nu$ of type $(2,0)$. Then with $\psi=\theta+\nu$ we want 
$$(d^*d(\rho+\psi +\bar\psi))^{2,0}=0$$
or, decomposing into types,
$$\bar\partial^*\partial\rho+\partial^*\partial \psi+\bar\partial^*\bar\partial\psi=0$$
or equivalently
\begin{equation}
(\partial^*\partial +\bar\partial^*\bar\partial)\psi=-\bar\partial^*\partial\rho
\label{oper}
\end{equation}
Now $\partial^*\partial +\bar\partial^*\bar\partial$ is elliptic. Indeed it
 is the sum of two non-negative second order operators $\partial^*\partial$ and 
 $\bar\partial^*\bar\partial$ and the latter is itself elliptic on $(2,0)$ forms.  
 The null space consists of $(2,0)$ forms $\psi$ satisfying $\bar\partial \psi=\partial\psi=0$ 
 (the holomorphic $2$-forms). It follows that  for each $\psi$ in this
 null-space 
$$\int_M(\bar\partial^*\partial\rho,\psi)=\int_M(\rho,\bar\partial\psi)=0$$
and so given $\rho$, we can solve the equation (\ref{oper}) for $\psi$. Using the Green's operator for $\partial^*\partial +\bar\partial^*\bar\partial$, we find $\theta \in L^2_{k+1}(\Lambda^2)$ as required.

Thus $\alpha_2=\alpha_1+d(\theta +\bar\theta) \in L^2_k(\Lambda^3)$ lies in the Banach space $E$, and the map $\alpha \mapsto \alpha_2$ is a continuous projection, orthogonal in $L^2$.
\vskip .25cm
If $\alpha$ is harmonic, then $d\alpha=0$ and $d^*\alpha=0$ and so in
particular 
$(d^*\alpha)^{2,0}=0$ and $\alpha\in E$. We  then split $E=E_1\oplus E_2$ where $E_1$ is the finite-dimensional space of harmonic forms and $E_2$ the exact ones in $E$.

 The Banach space $E$ is $L^2$-orthogonal to the orbit of $\Diff(M)$ through $\Omega$ and hence  is transversal to the orbits through a neighbourhood of $\Omega$. The functional $\Phi$ is smooth for $3$-forms $\alpha$ for which $\lambda(\alpha)<0$ at all points, so in order to define $\Phi$ on $E$ we need  uniform estimates on $\alpha$. The Sobolev embedding theorem tells us that in $6$ dimensions we can achieve this with $L^2_k(\Lambda^3)$ for $k>3$. Moreover in this range $L^2_k$ is a Banach algebra and so multiplication is smooth. Thus in a Sobolev neighbourhood of $\Omega$, $\Phi$ is a smooth function on $E$. Its derivative at $\alpha$ is given (from (\ref{firstv})) by
$$\delta \Phi(\dot\alpha)=-\int_M\hat\alpha\wedge \dot\alpha=\int_M(\ast\hat\alpha,\dot\alpha)\epsilon=\int_M(P(\ast\hat\alpha),\dot\alpha)\epsilon$$
where $P$ is orthogonal projection onto $E$.

The space $E$ is transverse to the $\Diff(M)$ orbits, and   $\Phi$ is constant on these,
so its derivative is determined by its derivative as a function on $E$. We want
the critical points of $\Phi$ restricted to a cohomology class so if $P_2$ denotes orthogonal projection onto the exact forms $E_2$, our critical points are the zeros of the function $F:E\rightarrow E_2$ defined by
$$F(\alpha)=P_2(\ast\hat\alpha).$$

\subsection{Invertibility of the derivative}

To apply the implicit function theorem, we want the derivative $D_2F:E_2\rightarrow E_2$ to be invertible, where, as we calculated in the previous section,
$$D_2F(\dot\alpha)=P_2(\ast J \dot\alpha).$$
The second variation calculation (\ref{var}) shows that this is an
injection for exact $\dot\alpha$.

We now prove surjectivity. Note that if $\beta$ is of type $(2,0)+(0,2)$, $P_2(\ast J d\beta)=0$, so surjectivity for $d\beta\in L_k^2(\Lambda^3)$ implies surjectivity on the transversal $E_2$. We use the following lemma (the proof is somewhat shorter in the K\"ahler case):

\begin{lem}
\label{surjlemma}
 If $\gamma\in L_{k+1}^2(\Lambda^{1,1})$ is a real form  which satisfies $(d^*d\gamma)^{2,0}=0$, then there exist a real form $\rho \in L_{k}^2(\Lambda^3)$ with $d^*\rho=0$ and a complex form $\sigma \in L_{k+1}^2(\Lambda^{2,2})$ such that
$$d\gamma=\rho+\partial^*\sigma+\bar \partial^*\bar\sigma.$$
\end{lem}
\begin{prf}
The condition $(d^*d\gamma)^{2,0}=0$ is equivalent to \begin{equation}
\bar\partial^*\partial \gamma=0
\label{condition}
\end{equation}
 Using the Green's operator for the $\bar\partial$-Laplacian, we write
\begin{equation}
\partial\gamma=H+\bar\partial G \bar\partial^*\partial\gamma+\bar\partial^* G \bar\partial\partial\gamma=H+\bar\partial^* G \bar\partial\partial\gamma
\label{Green}
\end{equation}
from (\ref{condition}). Here $H$ is the $\bar\partial$-harmonic component.

Now if $\bar\partial \theta=0$, applying the $\partial\bar\partial$-lemma to $d\theta$ we have $d\theta=\partial\bar\partial\nu$ and so $\theta-\bar\partial\nu$ is closed. Moreover, by imposing the condition $\bar\partial^*\nu=0$, if $\theta \in L^2_k$, $\nu\in L^2_{k+1}$. Similarly if $\bar\partial^* \theta=0$, there is a form $\nu$ such that $\psi=\theta-\bar\partial^*\nu$ is coclosed: $d^*\psi=0$. Applying this to the harmonic part $H$, which satisfies $\bar\partial^* H=0$, we have a coclosed form $\psi=H -\bar\partial^*\nu$ and so from (\ref{Green}) a $(2,2)$ form $\sigma=G\bar\partial\partial\gamma-\nu$ such that
$$\partial\gamma=\psi+\bar\partial^*\sigma.$$ Adding on the complex conjugate we obtain
$$d\gamma=\rho + \bar\partial^*\sigma+\partial^*\bar\sigma$$
as required.
\end{prf}

To continue with surjectivity, take $d\gamma \in E_2$. Since $P_2(d\gamma)=0$ if $\gamma$ is of type $(2,0)+(0,2)$, we may assume that $\gamma$ is of type $(1,1)$. From the lemma,  we can write
\begin{eqnarray*}
d\gamma&=&\rho + \bar\partial^*\sigma+\partial^*\bar\sigma\\
&=&\rho+(\bar\partial^*+\partial^*)u+i(\bar\partial^*-\partial^*)v\\
\end{eqnarray*}
where $\sigma=u+iv$. Now $v$ is of type $(2,2)$ so $i(\bar\partial^*-\partial^*)v=J(\bar\partial^*+\partial^*)v=Jd^*v$ which gives
$$Jd^*v=d\gamma -\rho-d^*u.$$
Since $\rho$ is coclosed, $d^*(\rho+d^*u)=0$ and thus the form $\rho+d^*u$ is orthogonal to all  exact forms. Hence
$P_2(Jd^*v)=P_2(d\gamma)$ or
$$P_2(\ast Jd\ast v)=P_2(d\gamma)$$ 
and we have surjectivity. From the open mapping theorem the derivative is now invertible. 

\subsection{The geometry of the moduli space}

We can now apply the implicit function theorem for Banach spaces to deduce that 
in a sufficiently small neighbourhood $U'$ of $\Omega$, the subspace $f^{-1}(0)\cap U'$ is diffeomorphic 
to a neighbourhood  of $H(\Omega)\in E_1$, the space of harmonic $3$-forms. Equivalently, 
the natural projection $p:E\rightarrow H^3(M,\R)$ identifies $f^{-1}(0)\cap U'$ with a neighbourhood $U$ 
of the cohomology class $[\Omega]$.

Taking $k>4$ for the Sobolev space $L^2_k(\Lambda^3)$, we have enough regularity to 
use the proof of the Newlander-Nirenberg theorem in \cite{NN} to deduce that the 
critical points of $\Phi$ on $U'\subset E$ define a family of complex structures with 
trivial canonical bundle on $M$. 
We now have what we need: a family of critical points parametrized by $U\subset
H^3(M,\R)$. \vskip .25cm
In this real approach to Calabi-Yau manifolds, the first thing we reach is a
moduli space which is 
an open set in the real vector space $H^3(M,\R)$. We  expect  to see a complex structure on this  (though the reader might  also appreciate   Kodaira's recollection in \cite{K} ``At first we did not take notice of the fact that the parameter appearing in the definition of a complex manifold is in general a complex one..."). We shall see next how a complex structure arises in our formalism.

\begin{prp}
\label{SK}
 The open set $U\subset H^3(M,\R)$ has the structure of a special pseudo-K\"ahler manifold of Hermitian signature $(1,h^{2,1}-1)$.
\end{prp}
\begin{prf}
The construction is essentially induced from that on $\Lambda^3 W^*$ of Section 3. In fact if $M$ is a $6$-torus the two are exactly the same.

We made use of  the symplectic form
$$\omega(\Omega_1,\Omega_2)\epsilon=\Omega_1\wedge\Omega_2$$
on $\Lambda^3 W^*$
and here we use its integrated form to define a flat symplectic structure on $H^3(M,\R)$: take closed forms $\Omega_1,\Omega_2$ with cohomology classes $[\Omega_1],[\Omega_2]$ and define
$$\omega([\Omega_1],[\Omega_2])=\int_M \Omega_1\wedge\Omega_2=([\Omega_1]\cup[\Omega_2])[M].$$
By Poincar\'e duality this is non-degenerate.

If $\Omega$ defines a complex structure then so does $\cos{\theta}\, \Omega+\sin{\theta}\,\hat \Omega$ so we may as well assume that the neighbourhood $U$ is invariant by the circle action
$$[\Omega]\mapsto \cos{\theta}\, [\Omega]+\sin{\theta}\,[\hat \Omega].$$
This  vector field preserves the symplectic form and has Hamiltonian function
$$\Psi=([\Omega]\cup[\hat \Omega])[M].$$
From (\ref{funcphi}) this is essentially the critical value of the functional $\Phi$ at the critical point $\Omega$.

The derivative of the map $[\Omega]\mapsto [\hat\Omega]$ on $U$ defines an almost complex structure since $\hat{\hat\Omega}=-\Omega$, and it is integrable just as in Proposition \ref{integprop}.

All that remains is to determine the signature. Here we use the fact that the $\partial\bar\partial$-lemma implies that the cohomology has a $(p,q)$ decomposition. This is the argument we used in proving surjectivity: each $\bar\partial$-closed form can be made closed by adding on a $\bar\partial$-exact form. The determination of the sign is then just as in Proposition \ref{openprop}.
\end{prf}
This special pseudo-K\"ahler structure (see \cite{Cand} for its origins) again has the property that 
the circle action is an isometry giving the quotient by the $\C^*$ action the structure 
of a projective special K\"ahler manifold. This quotient forgets 
the choice of holomorphic $3$-form and is a complex manifold of dimension $h^{2,1}$ 
which by our construction parametrizes a family of complex structures on $M$. Since from Proposition \ref{Jprop}, 
$J$ acts as $i$ on $H^{2,1}$, this is, as a complex manifold, the usual Kuranishi
moduli 
space.

\section{Manifolds with holonomy $G_2$}

\subsection{Linear algebra}

We have seen above how the consideration of a diffeomorphism-invariant functional on the space of $3$-forms 
on a $6$-manifold leads in a natural way to complex threefolds with trivial
canonical bundle and the special structure of their moduli space. We conclude now by addressing the analogous question in $7$ dimensions.

Let $W$ be a $7$-dimensional real vector space, and $\Lambda^3W^*$ the
$35$-dimensional vector space  of alternating multilinear $3$-forms on $W$.
 Take $\Omega\in \Lambda^3W^*$ and $v,w\in W$ and the
  interior products $\iota(v)\Omega, \iota(w)\Omega\in  \Lambda^2W^*$. 
  Then define $B_{\Omega}:W\otimes W \rightarrow \Lambda^7W^*$
by
$$B_{\Omega}(v,w)=-\frac{1}{6}\iota(v)\Omega\wedge \iota(w)\Omega\wedge \Omega
.$$ Since $2$-forms commute this is symmetric in $v,w$. 
It is a symmetric bilinear form on $W$ with values in the one-dimensional space $\Lambda^7W^*$, so it defines a linear map
$K_{\Omega}:W\rightarrow W^*\otimes \Lambda^7W^*$
and we then have
$$\det K_{\Omega}\in (\Lambda^7W^*)^9.$$
Since the exponent $9$ is odd, if $\det K_{\Omega}\ne 0$ then it defines this
way 
an orientation  on $W$ and, relative to this orientation we  choose the positive  root $$(\det K_{\Omega})^{1/9}\in \Lambda^7W^*.$$
We then  have  an inner product naturally defined by $\Omega$: 
$$g_{\Omega}(v,w)=B_{\Omega}(v,w)(\det K_{\Omega})^{-1/9}$$
 and $\det K_{\Omega}^{1/9}$ is its volume form.
 
Let $w_1,\dots,w_7$ be a basis of $W$ and $\theta_1,\dots,\theta_7$ the dual basis of $W^*$. For the $3$-form 
\begin{equation}
\varphi=(\theta_1\wedge\theta_2-\theta_3\wedge\theta_4)\wedge \theta_5 +(\theta_1\wedge\theta_3-\theta_4\wedge\theta_2)\wedge \theta_6
+(\theta_1\wedge\theta_4-\theta_2\wedge\theta_3)\wedge \theta_7 +\theta_4\theta_6\theta_7
\label{gform}
\end{equation}
an easy calculation gives the standard Euclidean metric
\begin{equation}
g_{\varphi}=\sum_1^7 \theta_i^2 
\label{detG}
\end{equation}
 The  stabilizer of $\varphi$ in $GL(7,\R)$  therefore preserves a metric and so is compact. 
 It is (see \cite{Sal})  the $14$-dimensional Lie group $G_2$: 
 the metric and $3$-form define the structure constants of the multiplication
 of the imaginary octonions. Since
$$\dim GL(7,\R)-\dim G_2=49-14=35=\dim \Lambda^3W^*$$
 the form $\varphi$ has an open orbit  $U\subset \Lambda^3W^*$. We call a form in this orbit {\it positive}.
 \begin{definition} Fix a basis vector $\epsilon \in \Lambda^7W^*$, and for a positive form $\Omega$  define  $\phi(\Omega) \in \R$ by $$\phi(\Omega)\epsilon=\det K_{\Omega}^{1/9}$$
 \end{definition}
 Note that since $K_{\Omega}$ is homogeneous of degree $3$ in $\Omega$, $\phi(\Omega)$  is homogeneous of degree $3\times 7/9=7/3$.

There is a formula for the invariant $\phi$ analogous to (\ref{funcphi}). Let $\ast:\Lambda^3W^*\rightarrow \Lambda^4W^*$ be the Hodge star operator of the metric $g_{\Omega}$, then
\begin{equation}
\Omega\wedge \ast \Omega=6\phi
\label{phistar}
\end{equation}

This follows from consideration of the standard $3$-form $\varphi$. As in \cite{Sal}, we have 
\begin{eqnarray*}
\ast\varphi&=&\theta_1\wedge\theta_2\wedge\theta_3\wedge\theta_4-(\theta_1\wedge\theta_2-\theta_3\wedge\theta_4)\wedge\theta_6\wedge\theta_7
-(\theta_1\wedge\theta_3-\theta_4\wedge\theta_2)\wedge\theta_7\wedge\theta_5\\
&-&(\theta_1\wedge\theta_4-\theta_2\wedge\theta_3)\wedge\theta_5\wedge\theta_6
\end{eqnarray*}
We want to find the derivative of $\phi:U\rightarrow \Lambda^3W^*$. The exterior product defines a non-degenerate pairing 
$$\Lambda^3W^*\otimes \Lambda^4W^*\rightarrow \Lambda^7W^*$$
and having fixed $\epsilon\in \Lambda^7W^*$ this defines an isomorphism
$$(\Lambda^3W^*)^*\cong \Lambda^4W^*.$$
Thus $d\phi$ takes values in $\Lambda^4W^*$. We have
  \begin{prp} 
  \label{718}
  The derivative of $\phi$ is given by:
$$d\phi(\dot\Omega)\epsilon=\frac{7}{18}\ast \Omega \wedge \dot\Omega$$
where $\ast: \Lambda^3 W^*\rightarrow \Lambda^4 W^*$ is the Hodge star operator of the metric $g_{\Omega}$.
\end{prp} 
\begin{prf}The group $G_2$ which stabilizes $\Omega$  has only one invariant in its action on  $\Lambda^3W^*$, namely the orthogonal projection onto $\Omega$ (see \cite{Sal}). Thus, since $\ast$ identifies $\Lambda^3W^*$ and $\Lambda^4W^*$ in an invariant way we can only have
$$d\phi(\dot\Omega)=C\ast \Omega \wedge \dot\Omega$$
for some real number $C$, possibly dependent on $\Omega$. 

But taking $\dot\Omega=\Omega$,   we have 
$$C(\ast \Omega \wedge \Omega)=d\phi(\Omega)=\frac{7}{3}\phi(\Omega)$$
since $\phi$ is homogeneous of degree $7/3$. From (\ref{phistar}), $C=7/18$.
\end{prf}

\subsection{The functional and its critical points}

We now define a functional in $7$ dimensions analogous to the one we have
analyzed in $6$ dimensions. Recall that if $\Omega$ is a $3$-form on a closed
 $7$-manifold $M$, then 
$\det K_{\Omega}^{1/9}$ is a section of $\Lambda^7T^*$. When it is
non-vanishing it is the volume form of a metric canonically determined by
$\Omega$. We define a functional $\Psi$ on $C^{\infty}(\Lambda^3T^*)$ by
 $$\Phi(\Omega)=\int_M \det K_{\Omega}^{1/9}.$$
We then have
\begin{thm} 
\label{G2thm}
Let $M$ be a closed $7$-manifold with a metric with holonomy $G_2$, with
defining $3$-form $\Omega$.
Then $\Omega$ is a critical point of the functional $\Phi$ restricted to the
cohomology class $[\Omega]\in H^3(M,\R)$.

 Conversely, if $\Omega$ is a critical point on a cohomology class of a closed 
 oriented $7$-manifold $M$ 
  such that $\Omega$ is everywhere positive, then $\Omega$ defines  on $M$ a
  metric with holonomy $G_2$.
  \end{thm} \begin{prf} Choose a non-vanishing $7$-form $\epsilon$ on $M$ and
  assume that $\Omega$ is a positive closed $3$-form.  Writing 
$$\det K_{\Omega}^{1/9}=\phi(\Omega)\epsilon$$
we have
$$\Phi(\Omega)=\int_M\phi(\Omega)\epsilon.$$
The first variation is:
$$\delta \Phi(\dot \Omega)=\int_Md\phi(\dot\Omega)\epsilon=\frac{7}{18}\int_M
\ast \Omega \wedge \dot\Omega$$
by (\ref{718}).  
Varying in a fixed cohomology class,  $\dot \Omega=d\varphi$ for some
$2$-form $\varphi$ so $\delta\Psi=0$ if and only if
$$\int_M d(\ast \Omega)\wedge \varphi=0$$
for all $\varphi$, or 
$$d(\ast \Omega)=0.$$
If $\Omega$ is covariant constant, so is $\ast \Omega$ and hence is closed.
Thus
a $G_2$-manifold gives a critical point of $\Psi$.

Conversely if $d\Omega=d(\ast\Omega)=0$ then $\Omega$ is covariant constant
for the metric $g_{\Omega}$ as shown by M.~Fern\'andez and A.~Gray  (see \cite
{Sal}).   \end{prf}
  \begin{ex}
  In D.~Joyce's foundational papers \cite{Joyce} there are many
  examples   of  compact $G_2$-manifolds. We can also take a $3$-dimensional
  Calabi-Yau manifold $N$ with $\Omega$ the real part of a holomorphic
  $3$-form   and $\omega$  the K\"ahler form and consider the product
  manifold $M=S^1\times N$. The $3$-form $\Omega +\omega\wedge d\theta$ defines
  a $G_2$-metric.
  \end{ex}
  
The second variation of $\Phi$ is the quadratic form 
$$\int_M D^2\phi(\dot\Omega_1,\dot\Omega_2)\epsilon$$
Now $18d\phi/7$ is the $4$-form $\Theta=\ast\Omega$ so for $D^2\phi$ what we need is the derivative of $\Theta(\Omega)$. This is provided
by the  lemma below (which is stated in \cite{Joyce}). To set it up, let us
consider the decomposition of $\Lambda^3W^*$ into $G_2$-representations. Since
$GL(7,\R)$ acts transitively on an open set of $\Lambda^3W^*$, we are looking
at the decomposition of the quotient $\lie{gl}(7,\R)/\lie{g}_2$ under the
action of the stabilizer $G_2$. Since $G_2\subset SO(7)$, and $G_2$ acts transitively on
$S^6$, the $21$-dimensional skew symmetric part of $\lie{gl}(7,\R)$ maps to the
$21-14=7$-dimensional representation in the quotient. The $28$-dimensional 
symmetric part decomposes into  multiples of the identity and the trace-free
symmetric tensors, which are $27$-dimensional. Denote as in \cite{Joyce} $\pi_1,\pi_7,\pi_{27}$
the projections onto these invariant subspaces, then we have:

\begin{lem}
\label{thetadiff}
 The derivative of the map $\Theta(\Omega)=\ast\Omega$ at $\Omega$ is
$$D\Theta(\dot\Omega)=\frac{4}{3}{\ast\pi_1(\dot\Omega)}+{\ast\pi_7(\dot
\Omega)}-{\ast\pi_{27}(\dot\Omega)}.$$ 
\end{lem}
\begin{prf} As in Proposition \ref{Jprop} we consider first the derivative applied to tangent vectors $X_a$ 
generated by $SL(7,\R)$, which span a subspace of codimension $1$. Since $\phi$
is invariant we have 
$$\sum_iX_a^i\frac{\partial \phi}{\partial x_i}=0$$
and so 
$$\sum_{i,j}\frac{\partial X^i_a}{\partial x_j}\frac{\partial \phi}{\partial
x_i}+X_a^i\frac{\partial^2 \phi}{\partial x_i\partial x_j}=0.$$
From (\ref{718}) and the fact that $DX_a=\rho(a)$ this can be written
$$D\Theta(X_a)=-\rho(a)\Theta$$
or
$$D\Theta(\rho(a)\Omega)=-\rho(a)(\ast\Omega).$$
From (\ref{action}), if we use an orthonormal basis $e_1,\dots,e_7$,
$$
\rho(a)(\ast\Omega)=\sum_{i,j}a_{ij}e_j\wedge\iota(e_i)\ast\Omega$$
Now use the relations $\iota(e_i)\alpha=\ast(e_i\wedge\ast\alpha)$ and $\ast^2=(-1)^{p(n-p)}$ on $p$-forms and we see that
 $$\sum_{i,j}a_{ij}e_j\wedge\iota(e_i)\ast\Omega=\ast\sum_{i,j}a_{ij}\iota(e_j)(e_i\wedge\Omega)=(\tr a) \ast\Omega-\ast\sum_{i,j}a_{ij}e_i\wedge\iota(e_j)\Omega.$$
Thus if $a_{ij}$ is skew-symmetric, $\rho(a)(\ast\Omega)=-\ast(\rho(a)\Omega)$
and so
\begin{equation}
D\Theta(\rho(a)\Omega)=-\rho(a)(\ast\Omega)
\label{skew}
\end{equation}
If $a_{ij}\in \lie{sl}(7,\R)$ it has trace zero, so if it is also symmetric  $\rho(a)(\ast\Omega)=\ast(\rho(a)\Omega)$ and 
\begin{equation}
D\Theta(\rho(a)\Omega)=\rho(a)(\ast\Omega)
\label{sym}
\end{equation}
There remains the unique direction orthogonal to the $SL(7,\R)$ orbit and this
is spanned by the Euler vector field $\Omega$. But
\begin{equation}
D\Theta(\Omega)_{\beta}=\frac{18}{7}\sum x_{\alpha}\frac{\partial^2
\phi}{\partial x_{\alpha}\partial x_{\beta}}=\frac{4}{3}\Theta_{\beta}=\frac{4}{3}\ast\Omega_{\beta}
\label{hom}
\end{equation}
since $\partial \phi/\partial x_i$ is homogeneous of degree $4/3$. From
(\ref{skew}),(\ref{sym}),(\ref{hom}) the lemma follows.
\end{prf}

\subsection{Nondegeneracy}

We  now prove, just as in the $6$-dimensional case, that at a $G_2$ manifold, the functional $\Phi$ has formally a Morse-Bott critical point. The argument is from one point of view like the proof of the $\partial\bar\partial$-lemma for K\"ahler manifolds, but is also in essence the proof of Proposition 3.3.2 in Joyce's paper \cite{Joyce}. 
\begin{prp} 
\label{Gnondegenprop}
Let $M$ be a compact $7$-manifold with a closed $3$-form $\Omega$ which defines a  $G_2$-structure. Then the Hessian of the functional $\Phi$ at $\Omega$ is nondegenerate transverse to the action of $\Diff(M)$ at $\Omega$.
\end{prp}
\begin{prf} The tangent space of the orbit of $\Diff(M)$ in the closed $3$-forms consists of the forms ${\cal L}_X\Omega=d(\iota(X)\Omega)$ where $X$ is a vector field. The map $X\mapsto \iota(X)\Omega$ embeds the tangent bundle in $\Lambda^2T^*$ as a  subbundle $\mu$ (in the notation of \cite{Sal}) and we have a $G_2$-invariant decomposition
$$ \Lambda^2T^*\cong \mu\oplus {\lie g}$$
where  ${\lie g}$ is the $14$-dimensional adjoint bundle of the $G_2$  structure. Thus the tangent space of the orbit consists of $3$-forms $d\chi$ where $\chi$ is a section of $\mu$. Defining $J:\Lambda^3T^*\rightarrow \Lambda^4T^*$ by  
\begin{equation}
J(\alpha)=\frac{4}{3}{\ast\pi_1(\alpha)}+{\ast\pi_7(\alpha)}-{\ast\pi_{27}(\alpha)}
\label{Jdefine}
\end{equation}
from (\ref{thetadiff}) we need to show that if $dJd\psi=0$ 
 then 
 $d\psi=d\chi$
 for $\chi\in C^{\infty}(\mu)$. By the diffeomorphism invariance of the functional any $\psi\in C^{\infty}(\mu)$ already satisfies $dJ(d\psi)=0$, so we need only restrict ourselves to $\psi \in C^{\infty}(\lie{g})$. We use the natural metric defined by the $G_2$ structure, so that nondegeneracy   in directions orthogonal to the $\Diff(M)$-orbits is the Morse-Bott condition. This means we must prove  that if 
$dJ(d\psi)=0$ and $\pi_{\mu}d^*d\psi=0$ then $d\psi=0$, where $\pi_{\mu}$ is the orthogonal projection onto the subbundle $\mu \subset \Lambda^2T^*$.

Assume we are given $\psi \in C^{\infty}(\lie{g})$   and $\pi_{\mu}d^*d\psi=0$. Using the Green's operator for the Hodge Laplacian we write
$$\psi=H(\psi)+dGd^*\psi+d^*Gd\psi$$
so that
$$d\psi=dd^*Gd\psi=dGd^*\psi$$
By assumption $\pi_{\mu}d^*d\psi=0$, so $d^*d\mu \in C^{\infty}(\lie{g})$. The Hodge Laplacian, and consequently also its Green's operator, commutes with orthogonal projections given by reduction of holonomy (see \cite{Besse} 1.141 for example). So $Gd^*d\psi \in C^{\infty}(\lie{g})$ and  $d\psi=d\chi$ where
\begin{equation}
\chi=Gd^*d\psi \in C^{\infty}(\lie{g})
\label{chi}
\end{equation}
and in particular $d^*\chi=0$. 

Now (as in \cite{Joyce}), since $\chi \in C^{\infty}(\lie{g})$ and $G_2$ acts as the irreducible $7$-dimensional representation on $\Lambda^6T^*$, we have $\chi\wedge \ast\Omega=0$. Thus $d\chi\wedge \ast\Omega=0\in C^{\infty}(\Lambda^7T^*)$ since $\ast\Omega$ is closed, and this means that \begin{equation}
\pi_1 d\chi=0
\label{pi1}
\end{equation}
Another algebraic property of $\chi \in C^{\infty}(\lie{g})$ (again see \cite{Joyce}) is 
$$\ast\chi=-\chi\wedge \Omega$$
so $d^*\chi=0$ implies $d\chi\wedge \Omega=0\in C^{\infty}(\Lambda^6T^*)$ since $\Omega$ is closed. This means that 
\begin{equation}
\pi_7 d\chi=0.
\label{pi7}
\end{equation}
 
 Now suppose $dJd\psi=dJd\chi=0$. Since $\pi_1d\chi=\pi_7d\chi=0$ we see from (\ref{Jdefine}) that $J(d\chi)=4\ast d\chi/3$ and so
 $$dJd\chi=-d\ast d\chi=0.$$
 But by Stokes' theorem $d*d\chi=0$ implies $d\chi=d\psi=0$, so we have nondegeneracy of the critical point.
\end{prf}

 \subsection{The moduli space}
 
 Given this nondegeneracy, we can proceed just as in Section 6 to define a local moduli space which is an open set in $H^3(M,\R)$. With $k>7/2$ we define the Sobolev space
$$E=\{\alpha \in L^2_k(\Lambda^3): d\alpha=0\quad{\rm and}\quad  
\pi_{\mu}d^*\alpha=0\}.$$
and proceed as before, using the function
$$F(\alpha)=P_2(\ast_0 \ast \alpha)$$
where $\ast_0$ is the $G_2$ metric at the critical point, and as before $\ast$ depends on $\alpha$. The derivative on exact forms at the critical point is then $P_2(\ast J d\psi)$ and the formal argument of nondegeneracy just proved shows that this is injective. 

For surjectivity of the derivative we need to prove that if $d\psi$, with $\psi \in C^{\infty}(\lie{g})$, satisfies $\pi_{\mu}d^*d\psi=0$ then  there exists $\nu$ such that 
$$d\psi=P_2\left(\frac{4}{3}{\pi_1(d\nu)}+{\pi_7(d\nu)}-{\pi_{27}(d\nu)}\right)$$
But $d\psi=\pi_1(d\psi)+\pi_7(d\psi)+\pi_{27}(d\psi)$ so from (\ref{pi1}) and (\ref{pi7}) 
$$d\psi=\pi_{27}(d\psi)$$
and we can take $\nu=-\psi$.

We can thus use the implicit function theorem again. The only other result
we need is the regularity of Einstein metrics (see \cite{Besse}, 5.26) to pass from a Sobolev space to a $C^{\infty}$ metric, and for this we can take $k>9/2$ in the Sobolev space. 

We have here produced  one use of the implicit function theorem to produce the moduli space, as was suggested in \cite{Joyce}.

\subsection{Geometry of the moduli space}

The natural differential geometry on the moduli space is again determined by the critical value of our functional $\Phi$. From (\ref{phistar}) we have the function 
$$\Phi([\Omega])=\frac{1}{6}\int_M \Omega \wedge \ast\Omega=\frac{1}{6}([\Omega]\cup [\ast\Omega])(M)$$
on the open set $U\subset H^3(M,\R)$ parametrizing the $G_2$-structures. We have, by analogy with (\ref{SK})
\begin{prp}
\label{SKG}
 Let $M$ be an irreducible $G_2$-manifold, and  $U\subset H^3(M,\R)$ a local moduli space. Then the Hessian of the function $\Phi([\Omega])$ defines a metric on $U$ of signature $(1,b_3-1)$.
\end{prp} 
\begin{prf} From (\ref{718}) the metric is a positive multiple of 
$$\int_M D\Theta(\dot\Omega)\wedge \dot\Omega$$
where $\dot \Omega$ is an infinitesimal deformation of a $G_2$ structure, and from (\ref{thetadiff}) this is 
$$\int_M
\left(\frac{4}{3}{\ast\pi_1(\dot\Omega)}+{\ast\pi_7(\dot\Omega)}-{\ast\pi_{27}(\dot\Omega)}\right)\wedge \dot\Omega$$
Both $\dot\Omega$ and $D\Theta(\dot\Omega)$ are closed, so we can replace $\dot\Omega$ by its harmonic representative $\alpha$. But $\pi_7\alpha$   defines  a harmonic $1$-form and since the Ricci tensor is zero for a $G_2$-manifold, any harmonic $1$-form is covariant constant by Bochner's theorem. Thus irreducibility implies that the harmonic $3$-forms decompose as a $1$-dimensional space of multiples of $\Omega$ and its $b_3-1$-dimensional orthogonal space which is the image of $\pi_{27}$. Then 
$$([\alpha],[\alpha])=\frac{4}{3}\int_M\Vert\pi_1\alpha\Vert^2-\int_M\Vert\pi_{27}\alpha\Vert^2$$
as required.
\end{prf}
\begin{rmk}
        Note that
        $\ast$ identifies $\Lambda^3W^*$ and $\Lambda^4 W^*$ as $G_2$
        representation spaces, and consequently the orbit of the $4$-form 
        $\ast\Omega$ under $GL(7,\R)$ is also open. We could therefore rederive  the characterization of $G_2$-manifolds through an invariant functional on $4$-forms instead of $3$-forms. One approach to this is symplectic:
        we   fix $\epsilon \in
        \Lambda^7W^*$, and identify  $\Lambda^4 W^*\cong (\Lambda^3W^*)^*$ via
        the natural pairing. Then the product
        $$\Lambda^3W^*\times\Lambda^4 W^*\cong T^*(\Lambda^3W^*)$$ 
        and Proposition \ref{718} tells us that the set of pairs
        $(\Omega,\ast\Omega)\in \Lambda^3W^*\times\Lambda^4 W^*$ is the graph
        of the derivative of the function $18\phi/7 $ on $\Lambda^3W^*$ and so
        is Lagrangian. Identifying $\Lambda^3W^*\times\Lambda^4 W^*$ the other way round with $T^*(\Lambda^4 W^*)$,
        this submanifold is also the graph of the derivative of a function
$\psi$ on $\Lambda^4 W^*$, which is the Legendre transform of $18\phi/7$. 
Since $\phi$ is homogeneous of degree $7/3$ in $x_i$, it follows  that
$\xi_i$ is homogeneous of degree $4/3$. Thus $x_i={\partial \psi}/{\partial
\xi_i}$ is homogeneous of degree $3/4$ and $\psi$ of degree $7/4$ in $\xi_i$.
In particular since $6\phi=\Omega\wedge \ast\Omega$, we have 
$$6\phi(\Omega)=\sum_i x_i {\xi}_i=\sum_i \xi_i \frac{\partial\psi}{\partial \xi_i}=\frac{7}{4}\psi(\ast\Omega)$$
and so
the functional $\Phi$ on the $3$-forms $\Omega$ can  be replaced by the functional
$$\Psi(\Theta)=\frac{7}{24}\int_M \psi(\Theta)$$
on the $4$-form $\Theta=\ast\Omega$.

The precise algebraic form of $\psi(\Theta)$ is as follows. For $\Theta \in
\Lambda^4W^*$,
given $v,w\in W$, we have $\iota(v)\iota(w)\Theta\in \Lambda^2 W^*$. The map
$v\wedge w\mapsto \iota(v)\iota(w)\Theta$ 
 defines a linear map
$$C_{\Theta}:\Lambda^2W\rightarrow \Lambda^2 W^*.$$
The determinant of $C_{\Theta}$ is a vector in $L^2$ where $L$ is the top exterior power of $\Lambda^2 W^*$.
The highest exterior power of $\Lambda^2 V$ for any vector space $V$ of dimension
$n$ 
is canonically isomorphic to $(\Lambda^n V)^{n-1}$ so in our case $\det
C_{\Theta}$ is a well-defined element of  
$$(\Lambda^7 W^*)^{12}$$
 and  $\det C_{\Theta}^{1/12}$ is a volume form, which will give an invariant
integrand. A calculation using the standard forms
$\varphi$ and $\ast\varphi$ shows that $\psi$ is in fact a nonzero multiple of
$\det C_{\Theta}^{1/12}$.
\end{rmk}

\end{document}